\documentclass[12pt]{amsart}

\usepackage{graphics}
\usepackage{amsmath}
\usepackage{amsfonts}
\usepackage{amssymb}
\usepackage{amscd}

\input xy
\xyoption{all}

\newtheorem{theorem}{\bf Theorem}
\newtheorem{lemma}[theorem]{\bf Lemma}

\newtheorem{proposition}[theorem]{\bf Proposition}
\newtheorem{corollary}[theorem]{\bf Corollary}
\newtheorem{definition}[theorem]{\bf Definition}
\newcommand{\Hom}{\operatorname{Hom}}

\title{The homotopy theory of equivalence relations}

\author{Finnur L\'arusson}
\address{School of Mathematical Sciences, University of Adelaide, Adelaide SA 5005, Australia.} 
\email{finnur.larusson@adelaide.edu.au}



\date{30 October 2006.  Minor changes 22 April 2007.}

\begin{document}

\begin{abstract}  We give a detailed exposition of the homotopy theory of equivalence relations, perhaps the simplest nontrivial example of a model structure.  
\end{abstract}

\maketitle

\tableofcontents

\section{Introduction}

\noindent
Abstract homotopy theory, also known as homotopical algebra, is an abstraction of the homotopy theory of topological spaces.  Its fundamental concept is the notion of a model category, or a model structure on a category, introduced by Quillen \cite{Quillen} in 1967.  Model structures have since appeared and been applied in a growing variety of mathematical areas.  Abstract homotopy theory is a rather technical subject.  It is usually quite difficult to verify the defining properties of a model category.  This note is a detailed treatment of what I believe to be the simplest nontrivial example of a model category: the category of equivalence relations.  As far as I know, such a treatment has not appeared in the literature before.  

Equivalence relations have a role to play in applied homotopy theory in situations where two objects are either equivalent or not and that is all there is to it.  Sometimes there is no natural notion of objects being equivalent in more than one way or of higher equivalences; sometimes one may simply wish to ignore such additional features.  There are applications in my own work where more sophisticated structures such as groupoids and simplicial sets are unnecessarily complicated and using equivalence relations is the natural way to go.  Such applications will typically involve embedding a geometric category in a category of diagrams or sheaves of equivalence relations.

The purpose of this note is to provide a detailed account of the basics as a foundation for applications.  Also, students of abstract homotopy theory, reading, say, \cite{DwyerSpalinski}, \cite{GoerssJardine}, or \cite{Hirschhorn}, may find it useful to see the central ideas of homotopical algebra presented in a concrete example with much lighter technical baggage than the standard example, simplicial sets, requires.

We start by defining our model structure on the category of equivalence relations, adopting the standard definition for groupoids, and verifying Quillen's axioms.  The homotopy category turns out to be the category of sets: the subject is homotopy theory at the level of $\pi_0$ and yet it is not trivial.  We prove that the model structure is proper and combinatorial, but not cellular, unlike the usual model structures on simplicial sets and topological spaces.  There is an enrichment of the category in itself, which interacts well with the model structure, as shown by a version of Quillen's Axiom SM7.  We introduce the pointwise fibration structure on the category of presheaves of equivalence relations over a small category and derive its basic properties.  Finally, we prove that colimits of diagrams of equivalence relations preserve acyclic maps, so homotopy colimits are just ordinary colimits, point out that limits do not, and give an explicit construction and a discussion of homotopy limits.  

We have tried to make the paper as self-contained as possible.  For definitions and other background not provided in detail, we refer the reader to \cite{GoerssJardine} and \cite{Hirschhorn}.  We have proved those of our results that are elementary---and most of them are---using only basic set theory and category theory: for aesthetic reasons, to make the proofs more accessible, and to lay bare the elementary nature of these results.  We have avoided unnecessarily sophisticated machinery that in some cases would have yielded shorter proofs.  The main exceptions are results about categories of diagrams in $\mathcal E$ rather than $\mathcal E$ itself: for these, deeper category theory and homotopy theory is needed, and precise references are provided.

\section{The model structure}

\noindent
An equivalence relation on a set is most simply viewed as a partition of the set.  A partitioned set can be considered both as a topological space that is a disjoint union, with the coproduct topology, of nonempty spaces with a trivial topology, and as a groupoid in which there is at most one isomorphism from one object to another.  We usually denote by $\sim_X$ or simply $\sim$ the equivalence relation corresponding to a partition of a set $X$.  The quotient set $X/\!\!\sim$ is sometimes denoted $\overline X$.  A morphism $f:X\to Y$ of partitioned sets is a map of sets that takes equivalent elements to equivalent elements.  In other words, the image by $f$ of a class in $X$ is contained in a class in $Y$.  If $X$ and $Y$ are viewed as topological spaces, this is precisely the definition of a continuous map; if they are viewed as groupoids, this is precisely the definition of a morphism of groupoids.  The induced map $\overline X\to\overline Y$ is denoted $\bar f$.  The category $\mathcal E$ of partitioned sets (more precisely, the category in which an object is a partition of a set and an arrow is a morphism of partitioned sets) is a full subcategory of the category of topological spaces and of the category of groupoids.

To provide a proper set-theoretic foundation for our study, we assume that a Grothendieck universe has been chosen and that all sets under consideration are elements of this universe (see \cite{Mac Lane}, Sec.\ I.6).  A category is called small if its objects form a set.

\begin{definition}  A map $f:X\to Y$ of partitioned sets is said to be:
\begin{enumerate}
\item  a cofibration if $f$ is injective;
\item  a fibration if $f$ maps each class of $X$ onto a class of $Y$;
\item  a weak equivalence if $f$ induces a bijection of quotient sets.
\end{enumerate}
\label{threeclasses}
\end{definition}

A weak equivalence is also called an equivalence or an acyclic map.  These notions agree with the usual ones for groupoids, introduced by Anderson in \cite{Anderson}; for more details, see \cite{Hollander} or \cite{Strickland}.  From the topological point of view, a fibration of partitioned sets as defined above is nothing but a Hurewicz fibration or equivalently a Serre fibration, and a weak equivalence is nothing but a topological weak equivalence or equivalently a homotopy equivalence.  A cofibration of partitioned sets, however, is more general than the two topological notions (both of which require a closed image, for example).

Here are three important maps of partitioned sets: cofibrations $i_0$ and $i_1$, and an acyclic cofibration $j$.
$$\xymatrix{ \varnothing \ar[r]^{i_0} & *+[F]{\cdot}  &  *+[F]{\boxdot\ \boxdot} \ar[r]^>>>>{i_1} &  *+[F]{\cdot\ \cdot} & *+[F]{\cdot} \ar[r]^{j}  &  *+[F]{\cdot \ \cdot}}$$

\medskip\noindent
The following result is immediate.  It almost says that the model structure we are about to define is cofibrantly generated (an additional set-theoretic regularity property is needed: see Theorem \ref{cofgenetc}).

\begin{proposition}  
\begin{enumerate}
\item  A map of partitioned sets is a fibration if and only if it has the right lifting property with respect to $j$.
\item  A map of partitioned sets is an acyclic fibration, that is, an acyclic surjection, if and only if it has the right lifting property with respect to $i_0$ and $i_1$.
\end{enumerate}
\label{rightlifting}
\end{proposition}

The next theorem is the main result of this section.

\begin{theorem}  The category $\mathcal E$ of partitioned sets with the three classes of maps defined above is a model category.  This means the following.
\begin{enumerate}
\item  $\mathcal E$ has all small limits and colimits.
\item  The two-out-of-three property:  If $f$ and $g$ are composable maps such that two of $f$, $g$, and $g\circ f$ are acyclic, then so is the third.
\item  If $f$ is a retract of $g$, and $g$ is acyclic, a cofibration, or a fibration, then so is $f$.
\item  Every commuting square
$$\xymatrix{A \ar[r]^{f} \ar[d]_{j} & X \ar[d]^{p} \\ B \ar[r]^{g} \ar@{-->}[ur]^{s} & Y}$$
where $j$ is a cofibration and $p$ is a fibration and one of them is acyclic, has a lifting $s$ making the two triangles commute.
\item  Every map can be functorially factored into a cofibration followed by an acyclic fibration, and into an acyclic cofibration followed by a fibration.
\end{enumerate}
\label{modelstructure}
\end{theorem}

\begin{proof}  (1)  The limit of a small diagram of partitioned sets is the set-limit with the coarsest partition (the largest equivalence relation) making the maps from the limit to the sets in the diagram morphisms of partitioned sets.  The colimit of a small diagram of partitioned sets is the set-colimit with the finest partition (the smallest equivalence relation) making the maps from the sets in the diagram to the colimit morphisms of partitioned sets.

(2)  Consider the maps $\bar f$ and $\bar g$ of quotient sets induced by $f$ and $g$.  Clearly, if two of $\bar f$, $\bar g$, and $\bar g\circ\bar f$ are bijections, then so is the third.

(3)  A retract of an injection is an injection.  If $f$ is a retract of $g$, then $\bar f$ is a retract of $\bar g$, and a retract of a bijection is a bijection.  As for fibrations, we observe that retractions preserve right lifting properties and invoke Proposition \ref{rightlifting}.

(4)  First, suppose $j$ is acyclic and let $b\in B$.  Since $j$ is acyclic, there is $a\in A$ with $b\sim j(a)$.  Then $g(b)\sim g(j(a))=p(f(a))$, so since $p$ is a fibration, $g(b)=p(x)$ for some $x\sim f(a)$.  If $b\in j(A)$, say $b=j(a)$ with $a\in A$, we take $x=f(a)$.  Set $s(b)=x$.  Then $p(s(b))=p(x)=g(b)$, so $s$ is a lifting in the square.  To verify that $s$ respects partitions, say $b\sim b'$ in $B$.  Find $a,a'\in A$ with $j(a)\sim b\sim b'\sim j(a')$.  Since $j$ is acyclic, $a\sim a'$, so $s(b)\sim f(a)\sim f(a')\sim s(b')$.

Next, suppose $p$ is acyclic, so $p$ is surjective, and let $b\in B$.  There is $x\in X$ with $g(b)=p(x)$.  If $b\in j(A)$, say $b=j(a)$ with $a\in A$, we take $x=f(a)$.  Set $s(b)=x$.  Then $p(s(b))=p(x)=g(b)$, so $s$ is a lifting in the square.  To verify that $s$ respects partitions, say $b\sim b'$ in $B$.  Then $p(x)=g(b)\sim g(b')=p(x')$ so $s(b)=x\sim x'=s(b')$ since $p$ is acyclic.

(5)  We imitate the constructions of mapping cylinders and mapping path spaces in topo\-logy.  Instead of the interval, we use the two-point set $I=\{0,1\}$ with $0\sim 1$.  Let $f:X\to Y$ be a map of partitioned sets.  Consider the diagram
$$\xymatrix{X \ar[r]^>>>>{\iota} \ar[d]_{f} & X\times I \ar[d] \ar[rdd]^{\tilde f} & \\ Y \ar[r] \ar@{=}[rrd] & M \ar[rd]|{p} & \\ & & Y}$$
in $\mathcal E$, where $M$ is the pushout $Y\cup_f(X\times I)$, $\iota(x)=(x,0)$, and $\tilde f(x,t)=f(x)$.  This gives a functorial factorization $f=p\circ j$, where $j:X\to M$, $x\mapsto (x,1)$.  Clearly, $j$ is injective, that is, a cofibration.  Also, $p$ is surjective, so to verify that $p$ is an acyclic fibration, we need to show that $\bar p$ is injective.  For $y\in Y$, $p^{-1}(y)$ is the image in $M$ of the subset $\{y\}\cup f^{-1}(y)\times I$ of $Y\amalg (X\times I)$.  This image lies in a single class in $M$, since if $x\in f^{-1}(y)$, then $y$ is identified in $M$ with $(x,0)$, which in turn is equivalent to $(x,1)$.  Hence, if $y'\sim y$, then $p^{-1}(y)$ and $p^{-1}(y')$ lie in a single class in $M$, so the preimage by $p$ of the class of $y$ is a single class in $M$.

Next, we define the path space $Y^I$ as the set of all maps $I\to Y$ of partitioned sets, that is, maps taking both $0$ and $1$ to the same class in $Y$, with two such maps considered equivalent if the corresponding classes are the same.  Note that this makes $Y$ and $Y^I$ weakly equivalent.  Consider the diagram
$$\xymatrix{X \ar@{=}[ddr] \ar[dr]|{i} \ar[drr]^{\tilde f} & & \\ & P \ar[d] \ar[r] & Y^I \ar[d]^{e} \\ & X \ar[r]^{f} & Y}$$
where $P$ is the pullback $X\times_f Y^I$, $e(\alpha)=\alpha(0)$, and $\tilde f(x)$ is the constant path taking both $0$ and $1$ to $f(x)$.  This gives a functorial factorization $f=q\circ i$, where $q:P\to Y$, $(x,\alpha)\mapsto \alpha(1)$.  Clearly, $i$ is injective.  To see that $i$ is acyclic, note that for $x,x'\in X$, we have $(x,\tilde f(x))=i(x)\sim i(x')=(x',\tilde f(x'))$ if and only if $x\sim x'$; also, if $(x,\alpha)\in P$, then $(x,\alpha)\sim(x,\tilde f(x))=i(x)$.  Finally, to verify that $q$ is a fibration, take $(x,\alpha)\in P$ (so $\alpha(0)=f(x)$) and let $y\sim q(x,\alpha)=\alpha(1)$.  Define $\beta:I\to Y$, $\beta(0)=f(x)$, $\beta(1)=y$.  Then $\beta(0)=f(x)=\alpha(0)\sim\alpha(1)\sim y=\beta(1)$, so $\beta\in Y^I$, $\beta\sim\alpha$, and $q(\beta)=y$.
\end{proof}

The initial object of $\mathcal E$ is the empty set $\varnothing$ (the empty colimit) and the final object of $\mathcal E$ is the one-point set $\ast$ (the empty limit), each with its unique partition.  Clearly, for every partitioned set $X$, the canonical map $\varnothing\to X$ is a cofibration and the canonical map $X\to\ast$ is a fibration, so $X$ is both cofibrant and fibrant, that is, bifibrant.  Thus, by the Whitehead Lemma (see \cite{Hirschhorn}, Thm.\ 7.5.10), every weak equivalence of partitioned sets is a homotopy equivalence.

There is a functor $Q$ from $\mathcal E$ to the category $\mathbf{Set}$ of sets, taking a partitioned set to its quotient set.  This functor has a right adjoint $R:\mathbf{Set}\to\mathcal E$, endowing a set with its discrete partition.  Namely, for a set $A$ and a partitioned set $X$, there is a natural bijection between maps $QX\to A$ and morphisms $X\to RA$.  (It is an exercise for the reader to show that $Q$ has no left adjoint.)  Now, $\mathbf{Set}$ has a rather trivial model structure in which the isomorphisms, that is, the bijections, are the weak equivalences and every map is both a cofibration and a fibration.  The pair $(Q,R)$ is then a Quillen pair, meaning that $Q$ preserves cofibrations and $R$ preserves fibrations (see \cite{Hirschhorn}, Sec.\ 8.5).  Furthermore, $(Q,R)$ is a pair of Quillen equivalences, meaning that the map $QX\to A$ is acyclic if and only if the corresponding morphism $X\to RA$ is acyclic.  This implies that the homotopy categories of $\mathcal E$ and of $\mathbf{Set}$ are equivalent; the latter is clearly $\mathbf{Set}$ itself, and we have proved the following result.

\begin{theorem}  The homotopy category of $\mathcal E$ is equivalent to the category of sets.
\label{htcat}
\end{theorem}

This does not mean that the homotopy theory of equivalence relations is trivial: there is more to a model structure than its homotopy category.  While a model structure is usually viewed as a tool for the study of the associated homotopy category, there are applications in which the model structure itself (in particular the fibrations and cofibrations) is the primary object of interest.  An example is the study \cite{Larusson} of lifting and extension properties in complex analysis.  Moreover, applications of the theory presented here will likely involve localizations of the pointwise fibration structure on categories of diagrams in $\mathcal E$ (see Corollary \ref{pointwisefibrationstructure}).  The homotopy categories of such localizations will typically be quite intricate.

\section{Good properties of the model structure}

\noindent
The first three theorems in this section show that the model structure on $\mathcal E$ has many of the good properties that we like model categories to have.  First, a model structure is said to be left proper if the pushout of an acyclic map along a cofibration is acyclic, right proper if the pullback of an acyclic map along a fibration is acyclic, and proper if it is both left proper and right proper.  

\begin{theorem}  The model category $\mathcal E$ is proper.
\label{proper}
\end{theorem}

\begin{proof}  Consider a commuting square
$$\xymatrix{A \ar[r]^{h} \ar[d]_{f} & B \ar[d]^{g} \\ C \ar[r]^{k} & D}$$
in $\mathcal E$.  First assume the square is a pushout, $h$ is a cofibration, and $f$ is acyclic.  We need to show that $g$ is acyclic.  We can take $D=(B\amalg C)/\!\!\approx$, where $\approx$ is the equivalence relation generated by letting $f(a)\approx h(a)$ for every $a\in A$.  The equivalence relation $\sim$ making $D$ a partitioned set is then generated by letting $g(b)\sim g(b')$ whenever $b\sim b'$ in $B$, and $k(c)\sim k(c')$ whenever $c\sim c'$ in $C$.

To show that $\bar g$ is injective, suppose $b, b'\in B$ and $g(b)\sim g(b')$.  By the descriptions of $\sim$ and $\approx$ just given, there is an even number $m\geq 0$ and points $b_1, \dots, b_m\in B$ and $c_1, \dots, c_m\in C$ such that 
$$b \sim b_1 \approx c_1 \sim c_2 \approx b_2 \sim b_3 \approx c_3 \sim \dots \sim c_m \approx b_m \sim b'.$$
Also, $b_j \approx c_j$ means that there is an odd number $i_j\geq 1$ and points $a_j^1,\dots, a_j^{i_j}\in A$ such that $b_j=h(a_j^1)$, 
$$f(a_j^1)=f(a_j^2), h(a_j^2)=h(a_j^3), f(a_j^3)=f(a_j^4), \dots, h(a_j^{i_j-1})=h(a_j^{i_j}),$$ 
and $c_j=f(a_j^{i_j})$.  Since $h$ is injective, $a_j^2=a_j^3$, $a_j^4=a_j^5, \dots, a_j^{i_j-1}=a_j^{i_j}$, so $f(a_j^1)=\dots=f(a_j^{i_j})=c_j$.  Thus, writing $a_j=a_j^1$, we now have 
$$b\sim h(a_1)\approx f(a_1) \sim f(a_2) \approx h(a_2) \sim h(a_3) \approx f(a_3) \sim \dots \sim f(a_m) \approx h(a_m) \sim b'.$$
Since $f$ is acyclic, we get $a_1\sim a_2$, $a_3\sim a_4, \dots, a_{m-1}\sim a_m$, so 
$$b\sim h(a_1) \sim h(a_2) \sim h(a_3) \sim \dots \sim  b'$$ 
and $b\sim b'$.

If $d\in D$ and $d\notin g(B)$, then $d=k(c)$ for some $c\in C$.  Since $f$ is acyclic, $c\sim f(a)$ for some $a\in A$.  Then $g(h(a))=k(f(a))\sim k(c)=d$.  This shows that $\bar g$ is surjective.

Right properness is easier.  Assume the square is a pullback, $k$ is a fibration, and $g$ is acyclic.  We need to show that $f$ is acyclic.  We can take $A=\{(b,c)\in B\times C:g(b)=k(c)\}$ with $(b,c)\sim(b',c')$ if and only if $b\sim b'$ and $c\sim c'$.  If $(b,c), (b',c')\in A$ and $c=f(b,c)\sim f(b',c')=c'$, then $g(b)=k(c)\sim k(c')=g(b')$, so $b\sim b'$ since $g$ is acyclic, and $(b,c)\sim(b',c')$.  Also, if $c\in C$, there is $b\in B$ with $g(b)\sim k(c)$ since $g$ is acyclic.  Since $k$ is a fibration, there is $c'\sim c$ with $k(c')= g(b)$; then $(b,c')\in A$ and $f(b,c')=c'\sim c$.
\end{proof}

If $X$ and $Y$ are partitioned sets, the set $\hom(X,Y)$ of morphisms $X\to Y$ carries an equivalence relation such that $f\sim g$ if $\bar f=\bar g$, that is, if $f(x)\sim g(x)$ for all $x\in X$.  We write $\Hom(X,Y)$ for the set $\hom(X,Y)$ with this equivalence relation.  Every composition map
$$\Hom(X,Y) \times \Hom(Y,Z) \to \Hom(X,Z), \quad (f,g)\mapsto g\circ f,$$
in particular the evaluation map
$$e:X\times\Hom(X,Y) \to Y, \quad (x,f)\mapsto f(x),$$
is a morphism of partitioned sets.  This enrichment of the category $\mathcal E$ in itself interacts well with the model structure on $\mathcal E$.  More precisely, we have the following version of Quillen's Axiom SM7.

\begin{theorem}   If $j:A\to B$ is a cofibration and $p:X\to Y$ is a fibration of partitioned sets, then the induced map
$$(j^*,p_*): \Hom(B,X) \to \Hom(A,X)\times_{\Hom(A,Y)} \Hom(B,Y)$$
is a fibration of partitioned sets, which is acyclic if $j$ or $p$ is acyclic.
\label{sm7}
\end{theorem}

We follow the usual method of proof for simplicial sets, as in \cite{GoerssJardine}, Sec.\ I.5.  The following lemma is called the Exponential Law.

\begin{lemma}  For partitioned sets $A$, $X$, and $Y$, the map 
$$e_*:\Hom(A,\Hom(X,Y)) \to \Hom (X\times A, Y), \quad e_*(g)(x,a)=g(a)(x),$$
is an isomorphism of partitioned sets, which is natural in $A$, $X$, and $Y$.
\label{exponentiallaw}  
\end{lemma}

\begin{proof}  The inverse morphism satisfies $e_*^{-1}(h)(a)(x)=h(a,x)$.
\end{proof}

\begin{proof}[Proof of Theorem \ref{sm7}]  Let $i:K\to L$ be a cofibration of partitioned sets.  By the Exponential Law, there is a lifting in a square of the form
$$\xymatrix{K \ar[r] \ar[d]_{i} & \Hom(B,X) \ar[d]^{(j^*,p_*)} \\
L \ar[r] & \Hom(A,X)\times_{\Hom(A,Y)} \Hom(B,Y)}$$
if and only if there is a lifting in the corresponding square
$$\xymatrix{(K\times B)\cup_{(K\times A)}(L\times A) \ar[r] \ar[d]_{\iota} & X \ar[d]^{p} \\
L\times B \ar[r] & Y}$$
Clearly, $\iota$ is injective and thus a cofibration.  To conclude the proof, we need to show that $\iota$ is acyclic if either $i$ or $j$ is.  Say $i$ is (the other case is analogous).  Then $i\times{\rm id}_B$ is acyclic and, by Theorem \ref{proper}, so is $({\rm id}_K\times j)_*(i\times{\rm id}_A)$.  The former map is $\iota$ precomposed by the latter map, so $\iota$ is acyclic by the two-out-of-three property.
\end{proof}

\begin{theorem}  The model category $\mathcal E$ is locally finitely presentable and cofibrantly generated, and hence combinatorial.
\label{cofgenetc}
\end{theorem}

\begin{proof}  Note that, first, every partitioned set is the colimit of the directed diagram of all its finite subsets, ordered by inclusion; second, a partitioned set $X$ is finitely presentable, meaning that $\hom(X,\cdot)$ preserves directed colimits in $\mathcal E$, if and only if it is finite; and, third, there is, up to isomorphism, only a set of finite partitioned sets.  Since $\mathcal E$ is cocomplete, this shows that $\mathcal E$ is locally finitely presentable (see \cite{AdamekRosicky}, Ch.\ 1).

The domains of the maps $i_0$, $i_1$, and $j$ in Proposition \ref{rightlifting} are finite and hence finitely presentable in $\mathcal E$, that is, $\aleph_0$-small relative to $\mathcal E$ in the language of \cite{Hirschhorn}.  Thus, $\mathcal E$ is cofibrantly generated with generating cofibrations $i_0$ and $i_1$ and a generating acyclic cofibration $j$ (see \cite{Hirschhorn}, Def.\ 11.1.2).  Finally, a locally presentable cofibrantly generated model category is, by definition, combinatorial.
\end{proof}

Every object in a locally presentable category is in fact presentable (small in the language of \cite{Hirschhorn}; see \cite{AdamekRosicky}, Prop.\ 1.16).  Hence, a locally presentable model category is cofibrantly generated, and thus combinatorial, if and only if its fibrations are characterized by a right lifting property with respect to a {\it set} of acyclic cofibrations and its acyclic fibrations are characterized by a right lifting property with respect to a {\it set} of cofibrations.  A good reference for basic facts on locally presentable model categories is \cite{Beke}, Sec.\ 1.  An important such fact is that the so-called {\it small object argument} works for any set of morphisms in a locally presentable category.  See also \cite{Dugger}, Sec.\ 2.

Let us say a few more words about the category theory of $\mathcal E$.  Let $\mathcal F$ be the small, full subcategory of $\mathcal E$ of finite partitioned sets.  The canonical functor from $\mathcal E$ to the category $\mathbf{Set}^{\mathcal F^{\text{op}}}$ of presheaves of sets on $\mathcal F$ is a full embedding, preserves directed colimits, and has a left adjoint, so it preserves limits (see \cite{AdamekRosicky}, Prop's 1.26, 1.27).  Hence, $\mathcal E$ is equivalent to a full, reflective subcategory of the topos $\mathbf{Set}^{\mathcal F^{\text{op}}}$, closed under directed colimits.  However, $\mathcal E$ itself is not a topos, if only because it lacks a subobject classifier.

We conclude this section with two results about the categories that one would actually use in geometric applications of the homotopy theory of equivalence relations.  The first result introduces the so-called pointwise fibration structure on categories of presheaves of equivalence relations.  A presheaf on a category $\mathcal C$ thought of as a site is the same thing as a diagram over the opposite category $\mathcal C^{\text{op}}$ thought of as an indexing category.  The two points of view are equivalent, but the former is more common in geometric applications.

\begin{corollary}  Let $\mathcal C$ be a small category and consider the category $\mathcal E^{\mathcal C^{\text{op}}}$ of presheaves of equi\-valence relations on $\mathcal C$.  There is a model structure on $\mathcal E^{\mathcal C^{\text{op}}}$ in which weak equivalences and fibrations are defined pointwise and cofibrations are defined by the left lifting property with respect to acyclic fibrations.  This model structure is proper and cofibrantly generated.  It is also locally finitely presentable and hence combinatorial.
\label{pointwisefibrationstructure}
\end{corollary}

\begin{proof}  The existence of the specified model structure and it being proper and cofibrantly generated follows from our previous theorems and the results of \cite{Hirschhorn}, Sec's 11.6, 13.1.  By \cite{AdamekRosicky}, Cor.\ 1.54, $\mathcal E^{\mathcal C^{\text{op}}}$ is locally finitely presentable since $\mathcal E$ is.
\end{proof}

The cofibrations in the pointwise fibration structure can be described somewhat explicitly: see \cite{Hirschhorn}, Thm.\ 11.6.1.

The category $\mathcal E^{\mathcal C^{\text{op}}}$ is enriched in $\mathcal E$ in much the same way that $\mathcal E$ itself is.  Namely, if $X$ and $Y$ are $\mathcal E$-valued presheaves on $\mathcal C$, the set $\hom(X,Y)$ of morphisms $X\to Y$ carries an equivalence relation such that $\phi\sim\psi$ if $\phi_C$ and $\psi_C$ are equivalent as maps $X(C)\to Y(C)$ of partitioned sets for every object $C$ in $\mathcal C$.  We write $\Hom(X,Y)$ for the set $\hom(X,Y)$ with this equivalence relation.  As before, every composition map
$$\Hom(X,Y) \times \Hom(Y,Z) \to \Hom(X,Z),$$
is a morphism of partitioned sets, and we have a version of Quillen's Axiom SM7.

\begin{theorem}  If $j:A\to B$ is a cofibration and $p:X\to Y$ is a fibration in $\mathcal E^{\mathcal C^{\text{op}}}$, then the induced map
$$(j^*,p_*): \Hom(B,X) \to \Hom(A,X)\times_{\Hom(A,Y)} \Hom(B,Y)$$
is a fibration in $\mathcal E$, which is acyclic if $j$ or $p$ is acyclic.
\label{presheafsm7}
\end{theorem}

To prove this, we start with a variant of the Exponential Law.  If $K$ is a partitioned set and $X$ is an object in $\mathcal E^{\mathcal C^{\text{op}}}$, we define the object $X^K$ in $\mathcal E^{\mathcal C^{\text{op}}}$ by setting $X^K(C)=\Hom(K, X(C))$ and letting a map $\rho:C\to D$ in $\mathcal C$ induce a map $\Hom(K, X(D)) \to \Hom(K, X(C))$ by postcomposition by $\rho^*:X(D)\to X(C)$.  This construction is easily verified to be covariant in $X$ and contravariant in $K$.

\begin{lemma}  For $\mathcal E$-valued presheaves $X$ and $Y$ on $\mathcal C$ and a partitioned set $K$, the map
$$\epsilon:\Hom(K,\Hom(X,Y))\to\Hom(X, Y^K), \quad \epsilon(g)_C(x)(k)=g(k)_C(x),$$
for each object $C$ of $\mathcal C$, $x\in X(C)$, and $k\in K$, is an isomorphism of partitioned sets, which is natural in $X$, $Y$, and $K$.
\label{presheafexponentiallaw}
\end{lemma}

\begin{proof}  The inverse morphism satisfies $\epsilon^{-1}(h)(k)_C(x)=h_C(x)(k)$.  \end{proof}

\begin{proof}[Proof of Theorem \ref{presheafsm7}]   Let $i:K\to L$ be a cofibration of partitioned sets.  By the Exponential Law, there is a lifting in a square of the form
$$\xymatrix{K \ar[r] \ar[d]_{i} & \Hom(B,X) \ar[d]^{(j^*,p_*)} \\
L \ar[r] & \Hom(A,X)\times_{\Hom(A,Y)} \Hom(B,Y)}$$
if and only if there is a lifting in the corresponding square
$$\xymatrix{ A\ar[r] \ar[d]_j & X^L \ar[d]^{q} \\
B \ar[r] & X^K\times_{Y^K} Y^L}$$
We need to verify that $q$ is a fibration which is acyclic if either $i$ or $p$ is.  Since the fibrations and weak equivalences in $\mathcal E^{\mathcal C^{\text{op}}}$ are defined pointwise, this follows directly from Theorem \ref{sm7}.
\end{proof}

\section{Effective monomorphisms}

\noindent
Cellularity is an important strengthening of cofibrant generation (see \cite{Hirschhorn}, Ch.\ 12).  The usual model structures on the category of simplicial sets and the category of topological spaces (say compactly generated and weakly Hausdorff) are both cellular.  We will show that $\mathcal E$ is not cellular.  (For another example of a cofibrantly generated model category that is not cellular, see \cite{Hirschhorn}, Ex.\ 12.1.7.)  Left proper cellular model categories admit left Bousfield localization with respect to any set of maps.  Fortunately, left proper combinatorial model categories do as well.

\begin{proposition}  $\mathcal E$ is not cellular.
\label{notcellular}
\end{proposition}

\begin{proof}  By the definition of cellularity, cofibrations in a cellular model category are effective monomorphisms (some say regular monomorphisms), that is, equalizers of pairs of maps (see \cite{Hirschhorn}, Def.\ 12.1.1).  This fails in $\mathcal E$.  For instance, the cofibration $i_1$ in Proposition \ref{rightlifting} is not an equalizer.  It is an easy exercise to show this directly.  It also follows from the fact that an effective monomorphism that is also an epimorphism is an isomorphism, whereas $i_1$ is both a monomorphism and an epimorphism, but not an isomorphism (the existence of such a map is another reason why $\mathcal E$ is not a topos).
\end{proof}

This result prompts us to take a closer look at effective monomorphisms in $\mathcal E$.  Note that the monomorphisms in $\mathcal E$ are precisely the injections, that is, the cofibrations.

\begin{proposition}  A monomorphism $m:A\to B$ in $\mathcal E$ is effective if and only if it induces an injection of quotient sets, that is, $m(a)\sim m(a')$ in $B$ if and only if $a\sim a'$ in $A$.
\end{proposition}

\begin{proof}  Since $m$ is an injection, it is the set-limit of the diagram $B\rightrightarrows B\cup_A B$ with the two natural inclusions, and $m$ is the $\mathcal E$-limit of this diagram if and only if $A$ carries the largest equivalence relation making $A\to B$ a morphism of partitioned sets (see the description of limits in $\mathcal E$ in the proof of Theorem \ref{modelstructure}).  This last condition is also implied by $m$ being the equalizer of any diagram $B\rightrightarrows C$, again by the description of limits in $\mathcal E$.
\end{proof}

It is now natural to ask the following question.  Is there a model structure on $\mathcal E$ (cellular, one would hope) in which the cofibrations are the effective monomorphisms and the weak equivalences are the same as before?  The answer is no: the acyclic cofibrations would still be the same (injections that induce bijections of quotient sets), so the whole model structure would be the same.

\section{Homotopy limits and colimits}

\noindent
Very simple examples, such as the pullback square
$$
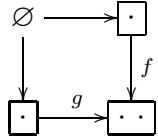
\xymatrix{\varnothing\ar[r] \ar[d] & *+[F]{\cdot} \ar[d]^{f} \\ *+[F]{\cdot} \ar[r]^{g} & *+[F]{\cdot\ \cdot} }$$
in $\mathcal E$, where $f$ and $g$ have different images, or the map of equalizers 
$$
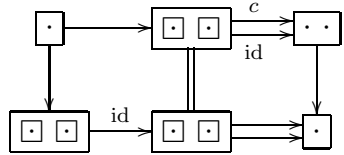
\xymatrix{*+[F]{\cdot} \ar[r] \ar[d] & *+[F]{\boxdot\ \boxdot} \ar@{=}[d] \ar@<-.5ex>[r]_{\rm id} \ar@<.5ex>[r]^{c} & *+[F]{\cdot\ \cdot} \ar[d] \\
*+[F]{\boxdot\ \boxdot} \ar[r]^{\rm id} & *+[F]{\boxdot\ \boxdot} \ar@<-.5ex>[r] \ar@<.5ex>[r] & *+[F]{\cdot}}$$
where $c$ is constant, show that limits of diagrams of partitioned sets need not preserve acyclicity: the map of limits induced by a pointwise acyclic map of diagrams need not be acyclic.  This forces us to study homotopy limits of diagrams of partitioned sets.  Indeed, one important aspect of model structures in general is that we can use the cofibrations and fibrations to construct and understand a modification of the ordinary notion of limits that does respect weak equivalences.  As for colimits, it is a special feature of the category of partitioned sets that they do preserve acyclicity, so homotopy colimits are just ordinary colimits.

\begin{theorem}  The morphism of colimits induced by a pointwise acyclic map of diagrams of partitioned sets is acyclic.
\label{colimits}
\end{theorem}

\begin{proof}  Every colimit is a coequalizer of a map of coproducts (see \cite{Mac Lane}, Sec.\ V.2).  Clearly, any set of weak equivalences of partitioned sets induces a weak equivalence from the coproduct of the sources to the coproduct of the targets.  Thus, we need to show that if we have a map of coequalizers of partitioned sets
$$\xymatrix{A \ar[d]_{\alpha}^{\sim} \ar@<-.5ex>[r]_{g} \ar@<.5ex>[r]^{f} & X \ar[d]_{\beta}^{\sim} \ar[r]^{p} & M \ar[d]_{\gamma} \\ B \ar@<-.5ex>[r]_{k} \ar@<.5ex>[r]^{h} & Y \ar[r]^{q} & N }$$
such that $h\alpha=\beta f$, $k\alpha=\beta g$, and $\alpha$ and $\beta$ are acyclic, then the induced map $\gamma$ is also acyclic.  We can take $M=X/\!\!\approx$, where $\approx$ is the smallest equivalence relation on $X$ with $f(a)\approx g(a)$ for all $a\in A$.  The equivalence relation $\sim$ on $M$ is the smallest one making $p$ a morphism, that is, the smallest one with $p(x)\sim p(x')$ if $x\sim x'$ in $X$.  Analogous remarks hold for $N$.

It is easy to show that $\bar\gamma$ is surjective: if $n\in N$, say $n=q(y)$, then, since $\bar\beta$ is surjective, there is $x\in X$ with $y\sim\beta(x)$, so $n=q(y)\sim q\beta(x)=\gamma p(x)$.

To show that $\bar\gamma$ is injective, let $m,m'\in M$ with $\gamma(m)\sim\gamma(m')$ in $N$.  Say $m=p(x)$, $m'=p(x')$, so $q\beta(x)\sim q\beta(x')$ in $N$.  Write $y=\beta(x)$, $y'=\beta(x')$.  The assumption that $q(y)\sim  q(y')$ in $N$ means that there is a string $y\sim y_1\approx y_1' \sim y_2 \approx\dots\approx y_\nu'\sim y'$ in $Y$.  The question is whether we can lift this string to a string joining $x$ and $x'$ in $X$.  Since $\bar\beta$ is injective, this is clear if $\nu=0$.  Assume $\nu\geq 1$ and consider the pair $y_1\approx y_1'$.  There are $b_1,\dots,b_j\in B$ and $y_1=z_0,z_1,\dots,z_j=y_1'\in Y$ such that $\{z_{i-1},z_i\}=\{h(b_i), k(b_i)\}$ for $i=1,\dots,j$.  One of the mutually analogous cases is when $y_1=h(b_1)$, $k(b_1)=h(b_2)$, \dots, $k(b_j)=y_1'$.  Since $\bar\alpha$ is surjective, there is $a_i\in A$ with $b_i\sim\alpha(a_i)$ for $i=1,\dots,j$.  Then $h(b_i)\sim h\alpha(a_i)=\beta f(a_i)$ and $k(b_i)\sim k\alpha(a_i)=\beta g(a_i)$, so $\beta(x)=y\sim y_1\sim\beta f(a_1)$, $\beta g(a_1)\sim k(b_1)=h(b_2)\sim \beta f(a_2)$, \dots, $\beta g(a_j)\sim k(b_j)=y_1'$.  Since $\bar\beta$ is injective, we get $x\sim f(a_1)$, $g(a_1)\sim f(a_2)$, \dots, $g(a_{j-1})\sim f(a_j)$.  If $\nu=1$, then $\beta g(a_j)\sim y_1'\sim y'=\beta(x')$, so $g(a_j)\sim x'$ and we have a string joining $x$ and $x'$, showing that $m\sim m'$.  If $\nu\geq 2$, we next consider the pair $y_2\approx y_2'$ and get $y_2\sim \beta f(\tilde a_1)$, say, so $\beta g(a_j)\sim y_1'\sim y_2\sim \beta f(\tilde a_1)$ and $g(a_j)\sim f(\tilde a_1)$, thus continuing the string that will eventually join $x$ and $x'$.
\end{proof}

Note that we did not need injectivity of $\bar\alpha$ in order to prove injectivity of $\bar\gamma$.  Indeed, if $A\neq\varnothing$, choose $s\in A$, adjoin a new element $t$ to $A$ such that $t\not\sim a$ for all $a\in A\setminus\{t\}$, and let $f$, $g$, and $\alpha$ take $s$ and $t$ to the same points in their respective targets.  Then the diagram is still a map of coequalizers, $\beta$ and $\gamma$ are still acyclic, but $\bar\alpha$ is not injective any more.

Now we turn to homotopy limits.  Every limit is an equalizer of a map of products (see \cite{Mac Lane}, Sec.\ V.2).  It is easy to verify that any set of weak equivalences of partitioned sets induces a weak equivalence from the product of the sources to the product of the targets.  We can therefore restrict our attention to equalizers.  The general theory of homotopy limits is quite involved and is developed in detail in \cite{Hirschhorn}, Ch's 18, 19.  See also \cite{DwyerSpalinski}, Sec.\ 10.  Our definition of homotopy equalizers in $\mathcal E$ is motivated by the general theory and justified by the results that follow.

\begin{definition}  The homotopy equalizer of a diagram of partitioned sets
$$\xymatrix{A \ar@<-.5ex>[r]_{g} \ar@<.5ex>[r]^{f} & X}$$
is the set $\{a\in A:f(a)\sim g(a)\}$ with the equivalence relation induced from $A$.
\label{homotopyequalizerdef}
\end{definition}

Let $\mathcal D$ be the indexing category $\xymatrix{ \bullet \ar@/_/[r] \ar@/^/[r] & \bullet }$.  Then the functor category $\mathcal E^{\mathcal D}$ is the category of diagrams $A\rightrightarrows X$ in $\mathcal E$ and maps between them.

\begin{theorem} 
\begin{enumerate} 
\item  The homotopy equalizer as defined above gives a functor $\mathcal E^{\mathcal D} \to \mathcal E$ with a natural transformation to the projection functor that takes $A\rightrightarrows X$ to $A$.
\item  The homotopy equalizer functor takes a pointwise acyclic map to an acyclic map.
\item  Let 
$$\xymatrix{A \ar@<-.5ex>[r]_{g} \ar@<.5ex>[r]^{f} & X}$$
be a diagram in $\mathcal E$.  The inclusion
$$\{a\in A:f(a)=g(a)\} \hookrightarrow \{a\in A:f(a)\sim g(a)\}$$
of the equalizer into the homotopy equalizer is acyclic if and only if for every $a\in A$ with $f(a)\sim g(a)$, there is $a'\in A$ with $f(a')=g(a')$ such that $a\sim a'$.
\end{enumerate}
\label{homotopyequalizer}
\end{theorem}

\begin{proof}  (1)  It is easily verified that a map in $\mathcal E^{\mathcal D}$ from $A\rightrightarrows X$ to $B\rightrightarrows Y$ yields a commuting diagram
$$\xymatrix{H \ar[r] \ar[d] & A \ar[d] \ar@<-.5ex>[r] \ar@<.5ex>[r] & X \ar[d] \\ K \ar[r] & B \ar@<-.5ex>[r] \ar@<.5ex>[r] & Y }$$
in $\mathcal E$, where we have denoted the homotopy equalizers of $A\rightrightarrows X$ and $B\rightrightarrows Y$ by $H$ and $K$ respectively, and $H\to A$ and $K\to B$ are the inclusions.

(2)  We need to show that if we have a map of homotopy equalizer diagrams
$$\xymatrix{H \ar[d]_{\gamma} \ar[r] & A \ar[d]_{\alpha}^{\sim} \ar@<-.5ex>[r]_{g} \ar@<.5ex>[r]^{f} & X \ar[d]_{\beta}^{\sim}  \\ K \ar[r] & B \ar@<-.5ex>[r]_{k} \ar@<.5ex>[r]^{h} & Y }$$
such that $h\alpha=\beta f$, $k\alpha=\beta g$, and $\alpha$ and $\beta$ are acyclic, then the induced map $\gamma$, obtained by restricting $\alpha$, is also acyclic.  Clearly, $\bar\gamma$ is injective since $\bar\alpha$ is.

To show that $\bar\gamma$ is surjective (this is what generally fails for ordinary equalizers), take $b\in B$ with $h(b)\sim k(b)$ and, using surjectivity of $\bar\alpha$, find $a\in A$ such that $\alpha(a)\sim b$.  Then
$$\beta f(a) = h\alpha(a) \sim h(b) \sim k(b) \sim k\alpha (a) = \beta g(a),$$
so $f(a)\sim g(a)$ since $\bar\beta$ is injective.  This shows that every element of $K$ is equivalent to an element in the image of $\gamma$.

(3)  The inclusion induces an injection of quotient sets by the definition of the equivalence relations on the equalizer and the homotopy equalizer: both are induced from $A$.  The given condition is precisely what it means for the inclusion to induce a surjection of quotient sets.
\end{proof}

There is a model structure on the diagram category $\mathcal E^{\mathcal D}$ in which the cofibrations and the weak equivalences are defined pointwise.  We call it the pointwise cofibration structure.  The relevance of such structures to homotopy limits is discussed in \cite{DwyerSpalinski}, Sec.\ 10.  (The pointwise fibration structure described in Corollary \ref{pointwisefibrationstructure} is similarly relevant to homotopy colimits.)  It is easy to verify that the characterization in Theorem \ref{homotopyequalizer} of when the homotopy equalizer of a diagram $A\rightrightarrows X$ is weakly equivalent to its ordinary equalizer means precisely that $A\rightrightarrows X$ (or rather the map from $A\rightrightarrows X$ to the final object $\ast\rightrightarrows 
\ast$ in $\mathcal E^{\mathcal D}$) has the right lifting property with respect to the map
$$\xymatrix{*+[F]{1} \ar[d]_{\alpha} \ar@<-.5ex>[r]_>>>>>>>{g} \ar@<.5ex>[r]^>>>>>>>{f} & *+[F]{1\quad 2} \ar[d]^{\beta} \\ *+[F]{1 \quad 2} \ar@<-.5ex>[r]_>>>>{k} \ar@<.5ex>[r]^>>>>{h} & *+[F]{1 \quad 2 \quad 3}}$$
where $f(1)=1$, $g(1)=2$, $\alpha(1)=1$, $\beta(1)=1$, $\beta(2)=2$, $h(1)=1$, $h(2)=3$, $k(1)=2$, $k(2)=3$, and the sources and targets have only one equivalence class each.  This map is an acyclic cofibration in $\mathcal E^{\mathcal D}$ with the pointwise cofibration structure, so the following corollary is immediate.

\begin{corollary}  If a diagram $A\rightrightarrows X$ is fibrant in the pointwise cofibration structure on $\mathcal E^{\mathcal D}$, then the natural map from its equalizer to its homotopy equalizer is acyclic.
\end{corollary}

Dually, for every small category $\mathcal C$ and every cofibrantly generated model category $\mathcal M$, the map of colimits induced by a pointwise acyclic map of diagrams in $\mathcal M^{\mathcal C}$ is acyclic if the diagrams are cofibrant in the pointwise fibration structure (see \cite{Hirschhorn}, Thm.\ 11.6.8).  This does not imply our Theorem \ref{colimits} because it is generally far from true that every diagram in $\mathcal E^{\mathcal C}$ is cofibrant in the pointwise fibration structure.  For example, with $\mathcal C=\mathcal D$ as above, the diagram $\xymatrix{*+[F]{\cdot}  \ar@<-.5ex>[r] \ar@<.5ex>[r]   & *+[F]{\cdot\ \cdot}}$\,, where the two arrows have the same image, has no map to the diagram $\xymatrix{*+[F]{\cdot}  \ar@<-.5ex>[r] \ar@<.5ex>[r]   & *+[F]{\cdot\ \cdot}}$\,, where the two arrows have different images, even though the map from the latter diagram to the final object $\xymatrix{*+[F]{\cdot}  \ar@<-.5ex>[r] \ar@<.5ex>[r]   & *+[F]{\cdot}}$ in $\mathcal E^{\mathcal D}$ is a pointwise acyclic fibration.

Let us consider the special case of homotopy pullbacks.  The pullback of a diagram 
$$\xymatrix{A \ar[r]^{f} & C & B \ar[l]_{g}}$$ 
is the equalizer of the diagram 
$$\xymatrix{A\times B \ar@<-.5ex>[r]_>>>>>{G} \ar@<.5ex>[r]^>>>>>{F} & C},$$ 
where $F(a,b)=f(a)$ and $G(a,b)=g(b)$.  One usually takes the homotopy pullback by replacing one of the maps, say $f$, with a fibration, meaning that one factors $f$ into an acyclic cofibration $A\to P$ followed by a fibration $P\to C$, and then taking the pullback $B\times_C P$.  Here, $P$ can be taken to be the mapping path space $A\times_f C^I$ described in the proof of Theorem \ref{modelstructure}.  This recipe for the homotopy pullback gives
$$B\times_C P = B\times_C (A\times_f C^I) = \{(b,a,\alpha)\in B\times A\times C^I: \alpha(0)=f(a), \alpha(1)=g(b)\}.$$
Recalling that an element of $C^I$, that is, a \lq\lq path\rq\rq\ in $C$, is simply a pair of equivalent elements of $C$, we see that this partitioned set is isomorphic to the set $\{(a,b)\in A\times B:f(a)\sim g(b)\}$, which is the homotopy equalizer of $A\times B\rightrightarrows C$ as we have defined it.

Finally, we remark that since limits and colimits in the presheaf category $\mathcal E^{\mathcal C^{\text{op}}}$ are taken pointwise, homotopy limits in the pointwise fibration structure of Corollary \ref{pointwisefibrationstructure} may be taken pointwise using Definition \ref{homotopyequalizerdef}.  Also, by Theorem \ref{colimits}, a pointwise acyclic map of diagrams of presheaves induces an acyclic map of their colimits.

\end{document}